\documentclass[oneside,12pt]{amsart}
\usepackage{amsmath,amsfonts, latexsym,amssymb}
\usepackage[active]{srcltx}

\theoremstyle{definition}

\numberwithin{equation}{section}


\begin{document}

\title[Variations of Lucas' Theorem 
Modulo Prime Powers]{Variations of Lucas' Theorem\\ Modulo Prime Powers} 

\author{Romeo Me\v strovi\' c}

\address{Maritime Faculty, University of Montenegro, Dobrota 36,
 85330 Kotor, Montenegro} \email{romeo@ac.me}

\subjclass {11A07, 11B65, 05A10.}
\keywords{Congruence modulo prime (prime power), 
Lucas' Theorem, a generalization of  Lucas' Theorem to prime powers, 
result of Jacobsthal, Wolstenholme prime.}

\maketitle

\begin{abstract}
Let $p$ be a prime, and let $k,n,m,n_0$ and $m_0$ be nonnegative integers 
such that $k\ge 1$, and $_0$ and $m_0$ are both less than $p$.
K. Davis and W. Webb  established that for a prime
$p\ge 5$ the following variation of Lucas' Theorem modulo prime powers
holds
 $$
{np^k +n_0 \choose mp^k+m_0}\equiv{np^{\lfloor(k-1)/3\rfloor}
\choose mp^{\lfloor(k-1)/3\rfloor}}
{n_0 \choose m_0}
\pmod{p^k}.
 $$
In the proof the authors used their earlier result
that present a generalized version of Lucas' Theorem.

In this paper  we present a  a simple
inductive  proof of the above congruence. Our proof is based 
on a classical congruence due to Jacobsthal, and we 
additionally use only some well known identities
for binomial coefficients.    
Moreover, we  prove that the assertion is also true for $p=2$ and $p=3$
if in the above congruence one replace
$\lfloor(k-1)/3\rfloor$ by $\lfloor k/2\rfloor$,
and by $\lfloor (k-1)/2\rfloor$, respectively.

As an application, in terms of  Lucas' type congruences, we 
obtain a new characterization of  Wolstenholme primes.

\end{abstract}
\vspace{6mm}

\centerline{1. INTRODUCTION AND MAIN RESULTS}

\vspace{2mm}

\par 
In 1878,  \'E.  Lucas proved a remarkable result which provides a 
simple way to compute
the binomial coefficient ${a\choose b}$ modulo a prime $p$ 
in terms of the binomial coefficients of the base-$p$ digits of 
nonnegative integers $a$ and $b$ with $b\le a$. Namely, 
if $p$ is a prime, and $n,m,n_0$ and $m_0$ are nonnegative integers
with $n_0,m_0\le p-1$, then a 
beautiful {\em theorem of Lucas} (\cite{l}; also see \cite{gr})
 states that for every prime $p$,
  $$
{np +n_0 \choose mp+m_0}\equiv {n \choose m}{n_0 \choose m_0}
\pmod{p}\eqno(1)
   $$
(with the usual convention that ${0\choose 0}=1$, and
${l \choose r}=0$ if $l<r$). After more than 110 years
D. F. Bailey established  that
under the above assumptions, $p$ can be replaced in (1) 
by $p^2$ \cite[Theorem 3]{b}, and  by $p^3$ 
if $p\ge 5$ \cite[Theorem 5]{b}. Moreover, it is noticed in \cite[p. 209]{b}
that in the congruence (1) $p$ cannot be replaced by $p^4$.
Using a Lucas' theorem for prime powers \cite[Theorem 2]{dw1}
(also cf. \cite[Theorem 2]{dw2}), in 1990 K. Davis and W. Webb 
\cite[Theorem 3]{dw2} generalized  Bailey's  congruences 
for any modulus $p^k$ with $p\ge 5$ and $k\ge 1$.
Their  result is improved quite recently
by the author of this paper in \cite{me}. 

Moreover, in 2007 Z.-W. Sun and D. M. Davis \cite{sd} and in 2009 
M. Chamberland and K. Dilcher \cite{chd} established analogues of 
Lucas' theorem for certain classes of binomial sums. 
Quite recently, the author of this article \cite{me4} discussed 
various cases of the congruences from Theorem A with $n_0=m_0=0$  
in dependence of different  values of exponents $k$ and $s$.

Another  generalization of  mentioned D. F. Bailey's Lucas-like theorem 
to every  prime powers $p^k$ with $p\ge 5$ and $k=2,3,\ldots$ was discovered 
in 1990 by K. S. Davis and W. A. Webb 
(\cite[Theorem 3]{dw1}, 
also see  \cite[p. 88, Theorem 5.1.2]{lo}) and independently 
by A. Granville \cite{gr1} (also see \cite{gr2} and 
 \cite[Theorem 1]{gr}).
Using mentioned result,  in 1993 K. S. Davis and W. A. Webb \cite{dw2}
generalized  Bailey's  congruences for any modulus 
$p^k$ with $p\ge 5$ and $k\ge 1$. Namely, they proved
the following congruence.\\

T{\scriptsize HEOREM} A (\cite[Theorem 3]{dw2}).
 {\em
 Let $p$ be any prime, and let $k,n,m,n_0$ and $s$ 
be positive  integers such that $0<n_0,m_0<p^s$. Then}
   $$
{np^{k+s} +n_0 \choose mp^{k+s}+m_0}\equiv {np^k\choose mp^k}
{n_0 \choose m_0}\pmod{p^{k+1}}.
    $$

R{\scriptsize EMARK} 1. As noticed above, Theorem A is proved 
by the authors  using their result in \cite[Theorem 3]{dw1}
which is slightly more complicated (cf. remarks by A. Granville
in \cite[Introduction]{gr}). The aim of this note is 
to give a simple elementary approach to the proof 
of Theorem A. For this purpose, in this note, we 
establish a simple induction proof 
 of Corollary of Theorem A (\cite[Corollary 1]{dw2}).
We point out that, proceeding by induction on $s$, 
the congruence in this Corollary (our Theorem given below)
allows us to establish a short and simple proof 
of Theorem A. This proof will be presented 
in the following version of this article.\\

T{\scriptsize HEOREM} 
(\cite[Corollary 1]{dw2}). {\em Let $p$ be any prime, 
and let $k,n,m,n_0$ and $m_0$  be nonnegative integers such that
$k\ge 1$, and $n_0$ and $m_0$ are both less than $p$.
If $p\ge 5$ then
 $$
{np^k +n_0 \choose mp^k+m_0}\equiv {np^{\lfloor(k-1)/3\rfloor}
\choose mp^{\lfloor(k-1)/3\rfloor}}
{n_0 \choose m_0}
\pmod{p^k},\eqno(2)
$$
where $\lfloor a\rfloor$ is  the greatest   integer less than or equal 
to $a$.\\\indent
Furthermore, for $p=2$ the congruence $(2)$ 
with $\lfloor k/2\rfloor$ instead of $\lfloor(k-1)/3\rfloor$
is satisfied, and for $p=3$ the congruence $(2)$ 
with $\lfloor (k-1)/2\rfloor$ instead of $\lfloor(k-1)/3\rfloor$
is also satisfied.}\\

As noticed above, the  congruences (2) for $k=2$ and 
$k=3$ are given by Bailey in \cite[Theorem 3 and Theorem 5, respectively]{b}
(our Corollaries  1 and 2, respectively).
Recall that  proof of Theorem 5 in \cite{b} 
is derived  by using  the congruence 
${np \choose mp}\equiv {n \choose m}\,(\bmod{\, p^3})$ 
with $p\ge 5$ \cite[Theorem 4]{b}
and a counting technique of M. Hausner from  \cite{h}. 
This theorem is refined modulo $p^5$ by a recent result of J. Zhao 
\cite[Theorem 3.5]{z}.\\\indent

 Our proof of the above theorem 
is inductive, and it is based on some congruences of Jacobsthal 
(see, e.g., \cite{gr}) and Sun and Davis \cite{sd}.
Namely,  the following lemma provides a basis for induction proof 
of Theorem.\\

 L{\scriptsize EMMA.}  {\em 
Let $n,m$ and $k$  be nonnegative integers with $m\le n$ and $k\ge 1$.
If $p$ is a prime greater than $3$, then
 $$
{np^k\choose mp^k}\equiv {np^{\lfloor(k-1)/3\rfloor} \choose 
mp^{\lfloor(k-1)/3\rfloor}}\pmod{p^k}.\eqno(3)
$$
Furthermore, for $p=2$ and $p=3$ we have
 $$
{n\cdot 2^k\choose m\cdot 2^k}\equiv {n\cdot 2^{\lfloor k/2\rfloor} \choose 
m\cdot 2^{\lfloor k/2 \rfloor}}\pmod{2^k},\eqno(4)
$$
$$
{n\cdot 3^k\choose m\cdot 3^k}\equiv {n\cdot 3^{\lfloor (k-1)/2\rfloor} 
\choose m\cdot 3^{\lfloor (k-1)/2 \rfloor}}\pmod{3^k}.\eqno(5)
$$}
\vspace{1mm}

{\em Proof.} We first  suppose that $p\ge 5$. 
Then we claim that the congruence  
$$
{np^k\choose mp^k} \equiv {np^{k-i}\choose mp^{k-i}} \pmod{p^{3(k-i+1)}}
\eqno(6)
   $$
holds for all  nonnegative integers $n,m,k$ and 
$i$ such that $1\le i\le k$.
If we put $i=k-\lfloor(k-1)/3\rfloor$ in (6), then since
$3(k-i+1)=3\lfloor(k-1)/3\rfloor +3\ge 3(k-3)/3+3=k$, 
we immediately obtain (3) from our Lemma.\\\indent
To prove (6),  we use induction on $i\ge 1$.   
By a result of Jacobsthal (see, e.g., \cite{gr}), 
 $$
{np\choose mp} \equiv {n\choose m} \pmod{p^e},\eqno(7)
   $$
for any integers $n\ge m\ge 0$ and prime  $p\ge 5$, where $e$ is the
power of $p$ dividing $p^3nm(n-m)$ (this exponent $e$ can only be 
increased if $p$ divides  $B_{p-3}$, the $(p-3)$rd {\em Bernoulli number}). 
Therefore, the congruence (7) with $np^{k-1}$ and $mp^{k-1}$ instead of  
$n$ and $m$, respectively, is satisfied for the exponent $e=3+3(k-1)=3k$.
This is in fact the congruence (6) with $i=1$.\\\indent
Now suppose that (6) holds for some $i$ such that $1\le i\le k-1$.
Then by a result of Jacobsthal mentioned above,
the congruence (7) with $np^{k-(i+1)}$ and $mp^{k-(i+1)}$ instead of  
$n$ and $m$, respectively, is satisfied for the exponent 
$e=3+3(k-(i+1))=3(k-i)$. This, together with the induction hypothesis 
given by (6), yields
   $$
{np^{k}\choose mp^{k}} \equiv {np^{k-(i+1)}\choose mp^{k-(i+1)}}  
\pmod{p^{3(k-i)}},
  $$
as desired.
\\\indent
If $p=2$ then by \cite[Lemma 3.2, the congruence (3.3)]{sd}, we have
 $$
{2n\choose 2m}\equiv (-1)^m {n\choose m}\pmod{2^{2{\rm ord}_2(n)+1}},
 $$
where ${\rm ord}_2(n)$ is the largest power of 2 dividing $n$.\\\indent
Then by induction on $k\ge 1$, similarly as above,
 easily follows   the congruence (4).\\\indent
Finally, if $p=3$ then by \cite[Lemma 3.2, the congruence (3.2)]{sd}, we have
 $$
{3n\choose 3m}\equiv  {n\choose m}\pmod{3^{2{\rm ord}_3(n)+2}},
 $$
where ${\rm ord}_3(n)$ is the largest power of 3 dividing $n$.\\\indent
Then by induction on $k\ge 1$ easily follows   the congruence (5).\\\indent
This completes the induction proof.\qed\\

{\em Proof of Theorem.} First suppose that $p\ge 5$,
and that $k$ is any fixed positive integer.
In order to prove the congruence (2), we proceed by induction on the sum
$s:=n_0+m_0\ge 0$, where $0\le n_0,m_0\le p-1$, and
hence $0\le s\le 2p-2$. If $s=0$, that is $n_0=m_0=0$, then the congruence (2)
reduces to the congruence (3) of our Lemma. \\\indent
Now suppose that the congruence (2) is satisfied for all $n,m,n_0$ 
and $m_0$ such that $n_0+m_0=s$ for some $s$ with $0\le s\le 2p-3$.
Next assume that $n_0$ and $m_0$ are any nonnegative integers such that 
$n_0+m_0=s+1$. Then consider the cases: $n_0<m_0$, $n_0=m_0\ge 1$  and 
$n_0\ge m_0+1$.\\\indent
{\tt Case 1.} $n_0<m_0$. Then ${n_0\choose m_0}=0$, and hence
the right side of (2) is equal to 0. Using the identity
${l\choose r}=\frac{l-r+1}{r}{l\choose r-1}$, we find that
   $$
{np^k+n_0  \choose mp^k+m_0}=\frac{p^k(n-m)-(m_0-n_0-1)}{mp^k+m_0}
{np^k+n_0  \choose mp^k+(m_0-1)}.
  $$
If $n_0=m_0-1$ then since $1\le m_0\le p-1$, the first factor on the right 
hand side of the above equality is divisible
by $p^k$. If $n_0<m_0-1$  then since $n_0+(m_0-1)=s$, 
by the induction hypothesis,  we get
 $$
{np^k+n_0  \choose mp^k+(m_0-1)}\equiv 
{np^{\lfloor(k-1)/3\rfloor}
\choose mp^{\lfloor(k-1)/3\rfloor}}
{n_0 \choose m_0-1}=0\pmod{p^k}.
  $$
Hence, in both cases we obtain 
  $$
{np^k+n_0  \choose mp^k+m_0}\equiv 0 = {np^{\lfloor(k-1)/3\rfloor}
\choose mp^{\lfloor(k-1)/3\rfloor}}
{n_0 \choose m_0}\pmod{p^k},
 $$
as desired.\\\indent
{\tt Case 2.} $n_0=m_0\ge 1$. 
If $n_0=m_0\ge 1$, then by the identity
${l\choose r}=\frac{l-r+1}{r}{l\choose r-1}$, in view of
$1\le n_0\le p-1$ and $n_0+(m_0-1)=s$, the induction hypothesis gives  
    $$
   \begin{array}{ll}
\displaystyle {np^k +n_0 \choose mp^k+n_0}&=\displaystyle
\frac{p^k(n-m)+1}{mp^k+n_0}
{np^k+n_0  \choose mp^k+(n_0-1)}\\[2.3ex]
&\displaystyle\equiv
\frac{p^k(n-m)+1}{mp^k+n_0}
{np^{\lfloor(k-1)/3\rfloor}
\choose mp^{\lfloor(k-1)/3\rfloor}}
{n_0 \choose n_0-1}\pmod{p^k}\\[2.3ex]
&=\displaystyle n_0\cdot\frac{p^k(n-m)+1}{mp^k+n_0}
{np^{\lfloor(k-1)/3\rfloor}
\choose mp^{\lfloor(k-1)/3\rfloor}}
\pmod{p^k}.
  \end{array}
    $$
This congruence and the fact that $1\le n_0\le p-1$ imply 
   $$
   \begin{array}{ll}
 & \displaystyle {np^k +n_0 \choose mp^k+n_0}-
{np^{\lfloor(k-1)/3\rfloor}\choose mp^{\lfloor(k-1)/3\rfloor}}
{n_0 \choose n_0}\\[2.3ex]
&\displaystyle\equiv 
\left(n_0\cdot \frac{p^k(n-m)+1}{mp^k+n_0}-1\right)
{np^{\lfloor(k-1)/3\rfloor}\choose mp^{\lfloor(k-1)/3\rfloor}}
\pmod{p^k}\\[2.3ex]
&\displaystyle= p^k\cdot\frac{n_0(n-m)-m}{mp^k+n_0}
{np^{\lfloor(k-1)/3\rfloor}\choose mp^{\lfloor(k-1)/3\rfloor}}
\equiv 0 \pmod{p^k},
\end{array}
    $$
whence follows (2).\\\indent
{\tt Case 3.} $n_0\ge m_0+1$.
Then we proceed in a similar way as in Case 2.
Using the identity ${l\choose r}=\frac{l}{l-r}{l-1\choose r}$,
in view of $1\le n_0-m_0\le p-1$ and $(n_0-1)+m_0=s$, 
the induction hypothesis yields    
   $$
   \begin{array}{ll}
  \displaystyle
{np^k +n_0 \choose mp^k+m_0}&\displaystyle =\frac{np^k+n_0}{p^k(n-m)+n_0-m_0}
{np^k+(n_0-1)  \choose mp^k+m_0}\\[2.3ex]
&\displaystyle\equiv
\frac{np^k+n_0}{p^k(n-m)+n_0-m_0}
{np^{\lfloor(k-1)/3\rfloor}
\choose mp^{\lfloor(k-1)/3\rfloor}}
{n_0-1 \choose m_0}\pmod{p^k}\\[2.3ex]
&\displaystyle=\frac{np^k+n_0}{p^k(n-m)+n_0-m_0}
{np^{\lfloor(k-1)/3\rfloor}
\choose mp^{\lfloor(k-1)/3\rfloor}}
{n_0\choose m_0}\cdot\frac{n_0-m_0}{n_0}\,.
\end{array}
  $$
The above congruence and the facts that $1\le n_0\le p-1$ and 
$1\le n_0-m_0\le p-1$, yield  
     $$
   \begin{array}{ll}
  &\displaystyle
 \quad {np^k +n_0 \choose mp^k+m_0}-
{np^{\lfloor(k-1)/3\rfloor}\choose mp^{\lfloor(k-1)/3\rfloor}}
{n_0 \choose m_0}\\[2.3ex]
&\displaystyle\equiv 
\left(\frac{n_0-m_0}{n_0}\cdot \frac{np^k+n_0}{p^k(n-m)+n_0-m_0}-1\right)
{np^{\lfloor(k-1)/3\rfloor}\choose mp^{\lfloor(k-1)/3\rfloor}}
{n_0\choose m_0}\pmod{p^k}\\[2.3ex]
&\displaystyle= p^k\cdot\frac{mn_0-nm_0}{n_0(p^k(n-m)+n_0-m_0)}
{np^{\lfloor(k-1)/3\rfloor}\choose mp^{\lfloor(k-1)/3\rfloor}}
{n_0\choose m_0}\equiv 0 \pmod{p^k},
  \end{array}
  $$
and so, (2) is satisfied. \\\indent
This concludes the assertion for any prime $p\ge 5$.\\\indent
The assertions of Theorem for $p=2$ and $p=3$ can be obtained
 by using the same method as in the above induction proof for $p\ge 5$, and 
hence may be omitted. Recall that the bases of induction proofs
related to $p=2$ and $p=3$ are the congruences (4) and (5) of  Lemma,
respectively. \\\indent
This completes the induction proof of Theorem. \qed\\

We now obtain two immediate consequences of Theorem.\\

C{\scriptsize OROLLARY} 1 (\cite[Theorem 3]{b}).   
{\em If $p$ is a prime, $n,m,n_0$ and $m_0$  are
 nonnegative integers, and  $n_0$ and $m_0$ are both less than $p$, then}
  $$
{np^2 +n_0 \choose mp^2+m_0}\equiv {n\choose m}
{n_0 \choose m_0}
\pmod{p^2}.
  $$
\vspace{1mm}

{\em Proof.} First observe that  the above assertion 
for $p\ge 5$ is a particular case
of Theorem for $k=2$. \\\indent
If $p=3$  then taking $k=2$ in (5) of  Leemma, we obtain
  $$
{9n \choose 9m}\equiv {n \choose m}\pmod{9}.
  $$
If we assume that the above congruence 
is a base of induction,  then applying the same method as in the proof 
of Theorem for the case $p\ge 5$, we obtain 
  $$
{9n +n_0 \choose 9m+m_0}\equiv {n \choose m}
{n_0 \choose m_0}
\pmod{9}, 
 $$
for all $n,m,n_0$ and $m_0$ with $0\le n_0\le 2$ and $0\le m_0\le 2$.\\\indent
Analogously, using the same argument, if we prove that  
   $$
{4n \choose 4m}\equiv {n \choose m}\pmod{4}, \eqno(8)
  $$
then it follows  that 
  $$
{4n +n_0 \choose 4m+m_0}\equiv {n \choose m}
{n_0 \choose m_0}
\pmod{4},
 $$
for all $n,m,n_0$ and $m_0$ such that $0\le n_0\le 1$ and $0\le m_0\le 1$.
\\\indent
To prove (8), note that by (4) of Lemma, we have
${4n \choose 4m}\equiv {2n \choose 2m}\,(\bmod{\, 4})$, 
and thus (8) is equivalent to the congruence 
 $$
{2n \choose 2m}\equiv {n \choose m}\pmod{4}\eqno(9) 
 $$
By the last congruence in the Proof 
of Lemma 3.2 in \cite{sd}, we have 
  $$
{2n\choose 2m} \equiv (-1)^m {n\choose m}
-(-1)^m 2n^2{n-1\choose m-1}(\frac{3+(-1)^m}{2}) 
\pmod{2^{2{\rm ord}_2(n)+2}}.\eqno(10) 
   $$  
If $m$ is even, then the above congruence immediately yields (9) for 
all $n$. If $m$ is odd and $n$ is even, then by Lucas' Theorem,
${n \choose m}\equiv 0 \,(\bmod{\, 2})$,
and thus (10) implies that
 $$
   \begin{array}{ll}
\displaystyle
{2n\choose 2m} &\displaystyle\equiv -{n\choose m}
+2n^2{n-1\choose m-1}
\pmod{4}\\[2.3ex]
&\displaystyle\equiv-{n\choose m}\equiv{n\choose m}\pmod{4}.
\end{array}
  $$\\\indent
Finally, if $n$ and $m$ are both odd, then from the identity  
$m{n\choose m}=n{n-1\choose m-1}$ we see that the integers
${n\choose m}$ and ${n-1\choose m-1}$ have the same parity.
This fact implies that 
$2{n \choose m}\equiv 2{n-1 \choose m-1}  \,(\bmod{\, 4})$, 
which together with the fact that $n^2\equiv 1 \,(\bmod{\, 4})$,
by (10) yields
  $$
{2n \choose 2m}\equiv - {n \choose m}+2{n-1 \choose m-1}\equiv
{n \choose m}\pmod{4}.
   $$\\\indent
This completes the proof.\qed\\

 C{\scriptsize OROLLARY} 2 (\cite[Theorem 5]{b}).  {\em 
Let $p$ be a prime greater than $3$.
If $n,m,n_0$ and $m_0$  are nonnegative integers with $n_0$ and $m_0$ 
less than $p$, then}
  $$
{np^3 +n_0 \choose mp^3+m_0}\equiv {n\choose m}
{n_0 \choose m_0}
\pmod{p^3}.
  $$
\vspace{1mm}

{\em Proof.} Clearly, the above assertion is a particular case
of  Theorem for $k=3$ with a prime $p\ge 5$. \qed

\vspace{6mm}

\centerline{2. A CHARACTERIZATION OF WOLSTENHOLME PRIMES}

\vspace{2mm}

A prime $p$ is said to be {\em Wolstenholme prime} if it 
satisfies the congruence ${2p-1\choose p-1} \equiv 1 \, (\bmod{\, p^4})$,
or equivalently,
 $$
{2p\choose p} \equiv 2 \pmod{p^4}.\eqno(11)
 $$
The two known such primes are 16843 and 2124679, and 
McIntosh and Roettger reported in \cite{mr} that these
primes are only two  Wolstenholme primes less than $10^9$.
However, McIntosh in \cite{m} conjectured
that there are infinitely many  Wolstenholme primes
(also see \cite{me2} and \cite[Section 7]{me3}).  \\\indent
As an application of Theorem of Section 1, 
in terms of  Lucas' type congruences, 
we obtain the following characterization of  Wolstenholme primes.\\

P{\scriptsize ROPOSITION.} {\em The following statements about a prime
$p\ge 5$ are equivalent.}\\\indent
$\,\,\,(i)$ $p$ {\em is a Wolstenholme prime};\\\indent
$\,\,(ii)$ {\em for  all  nonnegative integers $n$ and $m$,} 
   $$
{np\choose mp} \equiv {n\choose m} \pmod{p^4};\eqno(12)
    $$
\indent $(iii)$ {\em  for all nonnegative integers $n,m,n_0$ and $m_0$  such
that $n_0$ and $m_0$ are less than $p$,}
   $$
{np^4 +n_0 \choose mp^4+m_0}\equiv {n\choose m}
{n_0 \choose m_0}\pmod{p^4}.\eqno(13)
    $$
\vspace{1mm}

{\em Proof.} $(i)\Rightarrow (ii)$.
By  a special case of Glaisher's congruence 
(\cite[p. 21]{gl};   also cf. \cite[Theorem 2]{m}), for each  prime $p\ge 5$,
  $$
{2p-1\choose p-1} \equiv 1-\frac{2}{3}p^3B_{p-3} \pmod{p^4},
 $$
where $B_{p-3}$ is the $(p-3)$rd  Bernoulli number.
This shows that a prime $p$ is a Wolstenholme prime if and only if 
$p$ divides  the numerator of $B_{p-3}$.
On the other hand, by a result of Jacobsthal 
mentioned in the proof of Lemma (after the congruence (7)), 
the congruence (12) is satisfied
for any integers $n\ge m\ge 0$ and prime  $p\ge 5$ 
only if $p$ divides  $B_{p-3}$. \\\indent
$(ii)\Rightarrow (iii)$. 
Note that for any prime $p\ge 5$ and $k=4$ the congruence (2) 
of  Theorem becomes 
  $$
{np^4 +n_0 \choose mp^4+m_0}\equiv {np\choose mp}{n_0 \choose m_0}
\pmod{p^4}.
  $$
If we suppose that (12) is satisfied
for all  nonnegative integers $n$ and $m$, then 
(12) and the above congruence immediately yield (13), as desired.\\\indent
$(iii)\Rightarrow (i)$. If we suppose that (13) holds,
then taking $n=2$, $m=1$, $n_0=m_0=0$ in (13), we  obtain
the congruence ${2p^4\choose p^4} \equiv 2 \, (\bmod{\, p^4})$. 
On the other hand, taking $n=2$, $m=1$, $k=4$ and $i=3$
in (6), we have ${2p^4\choose p^4} \equiv {2p\choose p} \, (\bmod{\, p^6})$.  
These two congruences immediately imply (11), and
thus $p$ is a Wolstenholme prime.\\\indent
This completes the proof.\qed\\
  
R{\scriptsize EMARK} 2.  Note that 
for any prime $p\ge 5$ and for every $k\in\{4,5,6\}$ 
the congruence (2) of  Theorem becomes 
 $$
{np^k +n_0 \choose mp^k+m_0}\equiv {np\choose mp}
{n_0 \choose m_0}
\pmod{p^k}.\eqno(14)
  $$    
Note  that the first factor on the right side of (14)
is equal to ${np\choose mp}$, and that  for $k=4$ 
it can be replaced in (14) by ${n\choose m}$ if and only if 
${np\choose mp} \equiv {n\choose m}\ (\bmod{\,p^4})$.
Therefore, according to our Proposition, this is the case if and only if
$p$ is a Wolstenholme prime. Similarly, for $k=5$, this factor
can be replaced in (14) by  ${n\choose m}$ if and only if
 $$
{np\choose mp} \equiv {n\choose m} \pmod{p^5}\eqno(15)
    $$
for all $n$ and $m$. \\\indent
By {\em Wolstenholme's theorem} (see, e.g., \cite[Theorem 1]{z}),
if $p$ is a prime greater than 3, then the numerator of the 
fraction 
$$
H(p-1):=1+\frac{1}{2}+\frac{1}{3}+\cdots+\frac{1}{p-1}
 $$ 
is divisible by $p^2$. Now we define $w_p<p^2$ to be the unique 
nonnegative  integer such that $w_p\equiv H(p-1)/p^2\,\,(\bmod{\, p^2})$.
It is well known (see e.g., \cite{gl}) that  
    $$
w_p\equiv -\frac{1}{3}B_{p-3}\pmod{p}.
  $$
Furthermore, by a recent result of J. Zhao \cite[the congruence (10) of 
Theorem 3.2]{z}, for given prime $p\ge 7$ 
the congruence (15) is satisfied for all $n$ and $m$
if and only if $w_p=0$. However, using the argument based on the prime number 
theorem, McIntosh  \cite[p. 387]{m} conjectured that 
no  prime satisfies
the congruence ${2p-1\choose p-1} \equiv 1\,\, (\bmod{\,p^5})$.
Since the previous congruence is is a particular  case of (15) for $n=2$ and
$m=1$, McIntosh's Conjecture suggests  the following.\\

C{\scriptsize ONJECTURE.} {\em The exponent $\lfloor(k-1)/3\rfloor$ 
in the congruence $(2)$ of  Theorem can only be decreased for $k=4$ 
when $p$ is a Wolstenholme prime.}\\

R{\scriptsize EMARK} 3.
 Given any prime prime $p$ and $k\ge 2$, setting 
$n=m=n_0=1$ and $m_0=0$ in (2) of Theorem, we obtain
  $$
{p^k+1 \choose p^k}=p^k+1\equiv 1\pmod{p^k}.
 $$
This, together with the trivial fact that 
$p^k+1\not\equiv 1\,(\bmod{\,p^{k+1}})$,
shows that the exponent $k$ of the modulus 
$(\bmod \, p^k)$ in the congruence $(2)$ of  Theorem 
cannot  be increased for none  $k$ and $p$.

\end{document}